\documentclass[12pt,leqno]{article}
\usepackage{amsmath,amscd}

\def\int{\mathbb{Z}}

\def\gg{\mathfrak g}

\def\OO{{\cal O}}
\def\Oc{{\cal O}}

\def\ad{{\rm ad}}
\def\cc{\mathfrak{c}}
\def\ll{{\mathfrak l}}
\def\hh{{\mathfrak h}}

\def\mm{{\mathfrak m}}
\def\zz{{\mathfrak z}}
\def\jj{\mathfrak{J}}
\def\NN{{\mathcal N}}
\def\UU{{\mathcal U}}

\def\ZZ{{\mathbb Z}}
\def\CC{{\mathbb C}}

\def\pf{\proof}
\def\reg{{\textrm{reg}}}
\def\id{{\rm id}}

\def\Lie{{\rm Lie}}
\def\seq{\sim_\textrm{se}}
\def\loc{\sim_\textrm{loc}}
\def\epf{\hfill$\Box$ \medskip}
\usepackage{tikz}
\usepgflibrary{shapes.geometric}
\usetikzlibrary{matrix,arrows,decorations.pathmorphing}
\usepackage{amsmath,amsthm, amscd, amssymb, amsfonts, mathrsfs}
\usepackage[all]{xy}
\usepackage[vcentermath]{youngtab}
\usepackage{enumitem}
\usepackage{times}
\usepackage{graphics}
\usepackage{amssymb}
\usepackage{amscd}
\usepackage{amsmath}
\usepackage{color}
\usepackage{hyperref}

\title{Local geometry of Jordan classes in semisimple algebraic groups} 
\newtheorem{theorem}{Theorem}[section]
\newtheorem{lemma}[theorem]{Lemma}
\newtheorem{corollary}[theorem]{Corollary}
\newtheorem{proposition}[theorem]{Proposition}
\newtheorem{definition}[theorem]{Definition}
\newtheorem{remark}[theorem]{Remark}

\author{Filippo Ambrosio, Giovanna Carnovale, Francesco Esposito\\
Dipartimento di Matematica ``Tullio Levi-Civita''\\
Torre Archimede - via Trieste 63 - 35121 Padova - Italy\\
ambrosio@math.unipd.it, carnoval@math.unipd.it,\\
 esposito@math.unipd.it}

\begin{document}
\maketitle
\begin{abstract}
We prove that the closure of every Jordan class $J$ in a semisimple simply connected complex algebraic group $G$ at a point $x$ with Jordan decomposition $x=rv$ is smoothly equivalent to  the union of closures of those Jordan classes in the centraliser of $r$ that are contained in $J$ and contain $x$ in their closure. For $x$ unipotent we also show that the closure of $J$  around $x$ is smoothly equivalent to the closure of a Jordan class  in ${\rm Lie}(G)$ around $\exp^{-1}x$.  For $G$ simple we apply these results in order to determine a (non-exhaustive)  list of smooth sheets in $G$, the complete list of regular Jordan classes whose closure is normal and Cohen-Macaulay, and to prove that all sheets and Lusztig strata in ${\rm SL}_n({\mathbb C})$ are smooth. 
\end{abstract}

\section{Introduction}
Jordan classes in a reductive group or Lie algebra are locally closed, smooth, irreducible $G$-stable subsets of elements having similar Jordan decomposition.  

In the Lie algebra case they are also known as  decomposition classes and were introduced in \cite{BK}  in order to describe and parametrise sheets for the adjoint action of  a semisimple group on its Lie algebra. Their geometry has been studied in \cite{bo,broer,richardson,DR}. Sheets and birational sheets (recently introduced in \cite{losev}) are  unions of Jordan classes: these objects have a role in the representation theory of finite $W$-algebras, in the $G$-module structure of  rings of regular functions on adjoint orbits and in the description of  primitive ideals in enveloping algebras. The group version of Jordan classes made its first appearance in the work of Lusztig on the generalised Springer correspondence \cite{lusztig-inventiones}: they provide the stratification with respect to which character sheaves are constructible. Some of their properties and their closures have been studied in  \cite{gio-espo} in order to describe sheets for the action of a reductive group on itself.  Sheets in the group, in turn, are the irreducible components of the parts of the partition introduced in \cite{lustrata} as fibers of a map involving Springer representations with trivial local system, \cite{gio-MR}. 

Even though these Lie algebra and group stratifications were introduced to deal with distinct problems, they present similarities and it is natural to expect that the geometry of Jordan classes in a group and of Jordan classes in a Lie algebra are related. An example of the expected connection is to be found in \cite{gio-espo-normal} where the local geometry of the categorical quotient of the closure $\overline{J}$ of  a Jordan class in the group  $G$ has been related to the local geometry of categorical quotients of closures of Jordan classes in Lie algebra centralisers of semisimple elements contained in $\overline{J}$. This way the problem of normality or smoothness of $\overline{J}/\!/G$ could be related to the analogous problem for semisimple Lie algebras, whose solution is to be found in \cite{richardson, broer, DR}.

The first goal of this paper is to extend this approach to  the study of closures of  Jordan classes  in $G$ semisimple and simply connected. We prove in Theorem \ref{thm:main} that the closure of  a Jordan class $J$ in $G$ around a point $g$  with Jordan decomposition $g=rv$ is smoothly equivalent to  a union of  closures of Jordan classes in the centraliser of $r$ around the unipotent element $v$. We show that, up to a shift by $r$, the Jordan classes occurring  in this union are those classes containing $rv$ in their closure and contained in $J$ and we parametrise them in terms of Lie theoretic data depending on $J$ and $rv$. This allows to reduce the local study around any element to a local study around a unipotent one.
Then we prove in Theorem \ref{thm:stratification} that the exponential map identifies the Jordan stratification induced on a neighbourhood of the nilpotent cone in ${\rm Lie}(G)$ with the Jordan stratification induced on a neighbourhood of the unipotent variety in $G$,  preserving closure orderings. Therefore any closure of a Jordan class in $G$ containing a unipotent element $u$ is smoothly equivalent in the neighbourhood of $u$ to the closure of a Jordan class in ${\rm Lie}(G)$ in the neighbourhood of the logarithm of $u$.  We believe that these two  equivalences could establish new connections between representation theoretic objects related either to Jordan classes in a group or in a Lie algebra. For the present, we provide a series of applications to the study of geometry of closures of certain Jordan classes and of sheets.

Combining our local analysis and a theorem in \cite{richardson} describing when the closure of a regular Jordan class in a Lie algebra is normal and Cohen-Macaulay, we  prove that the closure of a regular Jordan class $J$  in $G$ is normal and Cohen-Macaulay if and only if  $\overline{J}/\!/G$ is normal if and only if  $\overline{J}/\!/G$ is smooth, Theorem \ref{thm:equivalence}. Since the list of classes $J$ for which $\overline{J}/\!/G$ is normal is known  \cite{gio-espo-normal}, this gives  the list of normal and Cohen-Macaulay closures of regular Jordan classes in $G$, see Remark \ref{cor:lista2}.

Every sheet $S$  contains a dense Jordan class $J_S$ and we provide necessary and sufficient conditions for a sheet  to be smooth, in terms of the local geometry of the closure of $J_S$. We also show in  Theorem \ref{thm:sheet-smooth-cod-one} that if $G$ is simple simply connected and classical and $\overline{J_S}/\!/G$ is normal in codimension 1,  then $S$ is always smooth, so since the list of  classes such that $\overline{J_S}/\!/G$ is normal in codimension 1 is known, \cite{gio-espo-normal}, we have a list of smooth sheets  for $G$ simple, simply-connected and classical, see Remark \ref{cor:lista-lisce}. We also provide the  list of smooth sheets when  $\overline{J_S}/\!/G$ is normal in codimension 1 for $G$ exceptional and simple in  Corollary \ref{cor:lista-lisce-ecc}.

When $G={\rm SL}_n({\mathbb C})$  
the situation is much simpler and we can conclude that all sheets and all Lusztig strata are smooth (Proposition \ref{prop:sln}). The general case is more involved and there are examples of singular and non-normal strata, for instance those containing the subregular unipotent conjugacy class when the root system is  doubly-laced. 

\section{Notation and preliminary results}\label{sec:notation}

In this and the following section $G$ is a complex connected reductive algebraic group; later it will be necessary to add further requirements on $G$. 
We fix a maximal torus $T$ with associated root system $\Phi$ and Weyl group $W$. We fix also a base $\Delta$ of $\Phi$ and $X_{\gamma}$, for $\gamma\in\Phi$ will be a root subgroup of $G$. If $\Phi$ is irreducible
 $\widetilde\Delta$ will be the union of $\Delta$ with the opposite of the highest root in $\Phi$;  the numbering of simple roots will be as in \cite{bourbaki}. We set $\mathfrak{g}={\rm Lie}(G)$, $\mathfrak{h}={\rm Lie}(T)$.
By abuse of terminology we will call Levi subalgebras (subgroups, respectively) the Levi subalgebras  (subgroups, respectively) of some parabolic subalgebra (subgroup, respectively) of $\gg$ (of $G$, respectively). 
The connected centralizers of semisimple elements in $G$ are called pseudo-Levi subgroups. If $\Phi$ is irreducible such groups are precisely those conjugate to a group of the form $G_\Pi=\langle T, X_{\pm\alpha},~|~\alpha\in\Pi\rangle$ for some $\Pi\subset \widetilde\Delta$, \cite[Proposition 3]{sommers},\cite[5.5]{lusztig}. The Weyl group of $G_\Pi$ will be denoted by $W_\Pi$ and $\mathfrak{g}_{\Pi}={\rm Lie}(G_\Pi)$. 

We use the dot to denote the conjugacy action of $G$ on itself, i.e., $h\cdot g=hgh^{-1}$. We denote by $Ad\colon G \to \rm GL(\gg)$ the adjoint representation of the group $G$ on $\gg$.  If $g\in G$, $V\subset G$ and $x\in\gg$,  we set
\begin{align*}
&C_G(g):=\{h\in G~|~ h\cdot g=hgh^{-1}=g\},\\ 
&C_G(V):=\bigcap_{v\in V}C_G(v),\\
&C_G(x):=\{h\in G~|~ Ad(h)(x)=x\},\\
&\cc_\gg(g):=\{y\in\gg~|~Ad(g)(y)=y\}=\Lie(C_G(g)),\\
&\cc_\gg(x):=\{y\in\gg~|~[x,y]=0\}=\Lie(C_G(x)).
\end{align*}

The conjugacy class of $g$ in a subgroup $H\leq G$ will be denoted by $H \cdot g = \OO^H_g$. For the adjoint orbit of $x \in \gg$, we use the notations $Ad(G)(x) = \mathfrak{O}^G_x$. 
If clear from the context, indices or superscripts will be omitted.
For any algebraic group $H$, the identity component  will be denoted by $H^\circ$ and the center  by $Z(H)$. The center of a subalgebra $\mm$ of $\gg$ will be denoted by $\zz(\mm)$.

When we write $g=su\in G$ or $x=x_s+x_n\in\gg$ we implicitly assume that $su$ ($x_s+x_n$, respectively) is the Jordan decomposition of $g$ ($x$, respectively), with $s$ semisimple and $u$ unipotent ($x_s$ semisimple and $x_n$ nilpotent, respectively).  We consider the elements in $\zz(\gg)$ as semisimple, so the semisimple part of $z+x$ for $z\in\zz(\gg)$ and  $x\in\gg$ is $z+x_s$.
 
If $G=Z(G)^\circ[G,G]$, we will  write $[G,G]_{sc}$ for the simply connected cover of the semisimple group $[G,G]$ and $G_{sc}:=Z(G)^\circ\times [G,G]_{sc}$. Also, 
$\pi\colon G_{sc}\to G$ will be a central isogeny and we will indicate by $T_{sc}$ the maximal torus in $G_{sc}$ such that $\pi(T_{sc})=T$.

It is well-known that the exponential map is a $G$-equivariant analytic map inducing  a $G$-equivariant analytic isomorphism between 
the nilpotent cone $\NN\subset\gg$ and the unipotent variety $\UU\subset G$, see \cite[\S 6.20]{hu-cc}.
For convenience in the exposition, we shall denote by $\exp_{sc}\colon \gg\to G_{sc}$ the scalar multiplication by $2\pi i$ followed by the exponential map
and by $\exp$ the composition $\pi\circ\exp_{sc}\colon \gg\to G$ and we shall call these maps the exponential maps. 

If $d\in\ZZ_{\geq0}$ and $X$ is a $G$-variety, we denote by $X_{(d)}$ the locally closed subset of $X$ consisting of points in orbits of dimension $d$. Their irreducible components are called the {\em sheets} for the action of $G$ on $X$.  For $Y\subset X$, we shall denote by $Y^{\reg}$ the set of points in $Y$ contained in $X_{(d)}$ for $d$ maximum such that $Y\cap X_{(d)}\neq\emptyset$.  For any group $H$ acting on a set $X$, unless otherwise stated, the stabiliser of $x\in X$ will be indicated by $H_x$.

For a surjective morphism $p\colon X\to Y$, we will say that $U\subset X$ is $p$-saturated  if $U=p^{-1}p(U)$.  We will use this notion for $X$ an affine $H$-variety with $H$ reductive and $Y=X/\!/H={\rm Spec}({\mathbb C}[X]^H)$ and we will denote the projection by $\pi_{X}$. In this case, $U$ is $\pi_{X}$-saturated if it is $H$-stable and such that  if $H\cdot u\subset U$, then $H\cdot x\subset U$ for every orbit $H\cdot x$ satisfying $\overline{H\cdot x}\cap\overline{H\cdot u}\neq\emptyset$, \cite[\S I]{luna}. For the main properties of the categorical quotient we refer to \cite[Theorem 1.24]{Brion}. 
See also \cite[\S 7.13]{Jantzen} for more details on the case $X=\mathfrak{g}$  and \cite[Chapter 3]{hu-cc} for the case of $X=G$, with adjoint action in both cases.

If $X\subset Y$ are topological spaces, we will denote by $\overline{X}^Y$ the closure of $X$ in $Y$. If the ambient space is clear, we will omit the superscript $Y$. We recall that when $X$ and $Y$ are algebraic varieties, the analytic closure coincides with the Zariski closure, \cite[Proposition 7]{serre} and that if $X$ is an algebraic variety and $x\in X$, then $X$ is unibranch, normal, smooth or Cohen-Macaulay at $x$ if and only if the corresponding analytic variety is so, \cite[Expos\'e XII, Proposition 2.1(vi), Proposition 3.1 (vii)]{SGA}.

\medskip
We will also need the following definition, see \cite[1.7]{he}
\begin{definition}Two pointed varieties $(X,x)$ and $(Y,y)$ are said to be smoothly equivalent if there exists a pointed variety $(Z,z)$ and two smooth maps $\varphi\colon Z\to X$ and $\psi\colon Z\to Y$ such that $\varphi(z)=x$ and $\psi(z)=y$. 
\end{definition}
Smooth equivalence is an equivalence relation on pointed varieties and it preserves the properties of being unibranch, normal, Cohen-Macaulay or smooth. We shall denote it by $\seq$. By \cite[Remark 2.1]{KP}
if $X$ and $Y$ are varieties satisfying $\dim Y=\dim X+d$, then  $(X,x)\seq(Y,y)$  if and only if $(X\times{\mathbb A}^d,(x,0))$ and $(Y,y)$ are locally analytically isomorphic. So, if $d=0$, then $(X,x)\seq(Y,y)$ if and only if there is a local analytic isomorphism in a neighbourhood of $x$ mapping $x$ to $y$: in this case we will also write $(X,x)\loc(Y,y)$. 
 
\section{Jordan classes and sheets in $G$ and $\mathfrak{g}$}\label{sec:Jordan_prelim}
In this section we recall  the necessary notions of Jordan classes in $\gg$ and $G$. For more information about them the reader is referred to \cite{bo,BK,broer} for the Lie algebra case and \cite{lusztig-inventiones,gio-espo} for the group case. The basic idea to keep in mind is that Jordan classes are irreducible subsets consisting of elements that, up to conjugation, have semisimple parts with same connected centraliser $M$ and nilpotent or unipotent part lying in the same $M$-orbit.

The Jordan class in $\gg$ containing the element $x=x_s+x_n$ is given by 
\begin{align}\jj_{\gg}(x):=Ad(G)(\zz(\cc_\gg(x_s))^{\reg}+x_n).
\end{align}
In other words it consists of all elements whose centralisers are $G$-conjugate to $\cc_\gg(x)$ \cite[39.1.6]{TY}. Jordan classes in $\gg$ are parametrised by $G$-orbits of pairs $(\ll,\mathfrak{O})$ where $\ll$ is a Levi subalgebra  of $\gg$  and $\mathfrak{O}$ is a nilpotent class in $\ll$. For the above class $\jj_{\gg}(x)$ we have $\ll=\cc_\gg(x_s)$ and $\mathfrak{O}= \mathfrak{O}_{x_n}\subset \ll$ and we will also indicate $\jj_\gg(x)$ by $\jj_\gg(\ll,\mathfrak{O})$.

The closure of $\jj_{\gg}(x)$ and its regular part are unions of Jordan classes and can be described as unions of adjoint orbits as follows:
\begin{align}\label{eq:closure-lie}
\overline{\jj_{\gg}(x)}&=\bigcup_{y_s\in\zz(\cc_\gg(x_s))}Ad(G)(y_s+\overline{{\rm Ind}_{\cc_\gg(x_s)}^{\cc_\gg(y_s)}\mathfrak{O}_{x_n}^{C_G(x_s)}})\\
\label{eq:reg-closure-lie}\overline{\jj_{\gg}(x)}^{\reg}&=\bigcup_{y_s\in\zz(\cc_\gg(x_s))}Ad(G)(y_s+{\rm Ind}_{\cc_\gg(x_s)}^{\cc_\gg(y_s)}\mathfrak{O}_{x_n}^{C_G(x_s)})
\end{align}
where ${\rm Ind}_{\cc_\gg(x_s)}^{\cc_\gg(y_s)}$ indicates Lusztig-Spaltenstein's induction of nilpotent orbits, \cite{lusp,bo}. Hence, $\jj_{\gg}(x)$ is closed if and only if $\zz(\cc_\gg(x_s))^{\reg}=\zz(\cc_\gg(x_s))=\zz(\gg)$ and $\mathfrak{O}_{x_n}^{C_G(x_s)}$ is closed, i.e., if and only if  $\jj_\gg(x)=\zz(\gg)$. 
It also follows from the above formula that the closure of any Jordan class in $\gg$ contains $0$, hence nilpotent elements. 

A Jordan class  $\jj'$ contained in $\overline{\jj_{\gg}(x)}^{\reg}$ is closed therein  if and only if $\overline{\jj'}^{\reg}=\jj'$ and this is the case if and only if $\jj'$ is the sum of $\zz(\gg)$ with the unique nilpotent orbit in $\overline{\jj(x)}^{\reg}$. 
 
\medskip

The Jordan class in $G$ containing the element $g=su$ is given by 
\begin{align} J_G(g):=G \cdot ((Z(C_G(s)^\circ)^\circ s)^{\reg}u).
\end{align}
The definition simplifies slightly if $G$ is simply connected, because $C_G(s)^\circ=C_G(s)$ for any semisimple element $s$. However, taking the connected component $Z(C_G(s)^\circ)^\circ s$ instead of $Z(C_G(s)^\circ)$  is necessary to guarantee irreducibility of a Jordan class. 
 
These classes are parametrised by $G$-orbits in the set ${\mathcal G}$ of triples $(M,Z(M)^\circ r,\Oc)$ where $M$ is a pseudo-Levi subgroup  of $G$; $Z(M)^\circ r$ is a coset in $Z(M)/Z(M)^\circ$ satisfying $C_G(Z(M)^\circ r)^\circ=M$  and $\Oc$ is a unipotent class in $M$. For the above class $J_G(g)$ we can take the triple:  $M=C_G(s)^\circ$, $Z(M)^\circ r=Z(M)^\circ s$, and $\Oc=\Oc^M_u$. We will denote $J_G(g)$ by $J_G(M, Z(M)^\circ s,\Oc)$. By construction, Jordan classes in $G$ are stable by left multiplication by elements in $Z(G)^\circ$.

The closure of $J_{G}(g)$ and its regular part are unions of Jordan classes and can be described as unions of conjugacy classes as follows:
\begin{align}\label{eq:closure-group}
\overline{J_{G}(M, Z(M)^\circ s,\Oc^M_u)}&=\bigcup_{z\in Z(M)^\circ s}G\cdot(z\overline{{\rm Ind}_{M}^{C_G(z)^\circ}{\OO}^M_{u})}\\
\label{eq:reg-closure-group}
\overline{J_{G}(M, Z(M)^\circ s,\Oc_u^M)}^{\reg}&=\bigcup_{z\in Z(M)^\circ s}G\cdot(z{\rm Ind}_{M}^{C_G(z)^\circ}{\OO}_{u}^M)
\end{align}
where ${\rm Ind}_{M}^{C_G(z)^\circ}$ indicates Lusztig-Spaltenstein's induction  of unipotent conjugacy classes, \cite{lusp}\cite[Proposition 4.8]{gio-espo}. The Jordan class $J_G(g)$ is closed if and only if 
$Z(M)^\circ s=(Z(M)^\circ s)^{\reg}=Z(G)^\circ s$ and $\overline{{\OO}_{u}^{M}}={\OO}_{u}^{M}$. One can verify that this happens if and only if $u=1$ and 
$M /Z(G)^\circ$ is semisimple, i.e., if and only if $g=s$ is semisimple and {\em isolated} in the terminology of \cite{lusztig-inventiones}. A Jordan class $J'$ contained in $\overline{J_G(M, Z(M)^\circ s,\Oc^M_u)}^{\reg}$ is closed therein if and only if $\overline{J'}^{\reg}=J'$ and this is the case if and only if 
$J'=J_G(M',Z(M')^\circ r,\Oc_v)\subset \overline{J_G(M, Z(M)^\circ s,\Oc_u^M)}^{\reg}$ with $M'/Z(G)^\circ$ semisimple, i.e., the semisimple part of the elements in $J'$ are  isolated.

It is worthwhile to notice that, in contrast with the Lie algebra situation, not all closures of Jordan classes contain a unipotent conjugacy class, even up to a central element. In fact, 
$\overline{J_G(M,Z(M)^\circ s,\Oc)}\cap Z(G){\mathcal U}\neq\emptyset$ if and only if $\overline{J_G(M,Z(M)^\circ s,\Oc)}\cap Z(G)\neq\emptyset$ and the latter holds  if and only if $M$ is a Levi subgroup. Also, $\overline{J_G(M,Z(M)^\circ s,\Oc)}\cap{\mathcal U}\neq\emptyset$ if and only if $1\in \overline{J_G(M,Z(M)^\circ s,\Oc)}$ if and only if $M$ is a Levi subgroup and 
$Z(M)^\circ s=Z(M)^\circ$, see formula \eqref{eq:closure-group} and the proof of \cite[Proposition 5.6]{gio-espo}.

\medskip

Using our choice of maximal torus $T$ we can simplify the parametrisation of Jordan classes in $G$ by reducing the set of admissible triples and the symmetry group acting on it. Let ${\mathcal T}=\{(M, Z(M)^\circ s,\Oc_u^M)\in{\mathcal G}~|~T\subset M\}$. Observe that in this case $Z(M)^\circ s\subset T$, that $N_G(T)$ acts on $\mathcal{T}$ and that $T$ acts trivially, so $W$ acts on $\mathcal{T}$. 

\begin{proposition}\label{prop:param}
Jordan classes in $G$ are parametrised by elements in ${\mathcal T}/W$. 
\end{proposition}
\pf We need to show that ${\mathcal G}/G$ is in bijection with ${\mathcal T}/W$. First of all, since all semisimple classes in $G$ have a representative in $T$, any triple  in $\mathcal{G}$ is $G$-conjugate to a triple in $\mathcal{T}$. We show that two triples in $\mathcal{T}$ are $G$-conjugate if and only if they lie in the same $W$-orbit. One direction is immediate, as $N_G(T)\subset G$. Let $(M_1,Z(M_1)^\circ s_1,\Oc_1)$ and $(M_2,Z(M_2)^\circ s_2,\Oc_2)\in\mathcal{T}$ and assume 
\begin{align*}(M_2,Z(M_2)^\circ s_2,\Oc_2)&=g\cdot(M_1,Z(M_1)^\circ s_1,\Oc_1),&\textrm{ for some }g\in G.\end{align*} Since all maximal tori in $M_2$ are $M_2$-conjugate, there exists $m\in M_2$ such that $\dot{w}:=mg\in N_G(T)$, and 
$(M_2,Z(M_2)^\circ s_2,\Oc_2)=\dot{w}\cdot(M_1,Z(M_1)^\circ s_1,\Oc_1)$
\epf

Jordan classes  in $\gg$ and $G$ form a partition of their ambient variety  into finitely many locally closed, irreducible,  smooth  $G$-stable subsets \cite{broer,BK,lusztig-inventiones}.  If the ambient Lie algebra or group is clear, we will omit the subscript ${\gg}$ or $G$.

\medskip

The sheets for the action of $G$ on $\gg$ or $G$ are obtained as follows, \cite{BK,gio-espo}: every sheet $S$ in $\gg$ (in $G$, respectively) contains a unique dense Jordan class $\jj$ ($J$, respectively) and $S=\overline{\jj}^{\reg}$ ($S=\overline{J}^{\reg}$, respectively). A Jordan class $\jj(\ll,\mathfrak{O})$ ($J(M,Z(M)^\circ s, \Oc)$, respectively) is dense in a sheet if and only if $\mathfrak{O}$ is {\em rigid} in $\ll$ ($\Oc$ is {\em rigid} in $M$, respectively), i.e., it is not induced from an orbit (conjugacy class, respectively) in a proper Levi subalgebra (subgroup, respectively).

\section{Reduction to unipotent elements}\label{sec:to_uni}

In this section $G$ is semisimple and simply connected. We begin our local study of Jordan classes. We will use a variant of Luna's \'etale slice theorem to reduce the study of the closure of a Jordan class in $G$ in the neighbourhood of an element $rv$ to the study of the closures of several Jordan classes in $C_G(r)$ in the neighbourhood of the unipotent part $v$.

We recall that if $H$ is a  reductive subgroup of $G$ acting on a variety $X$  then $G\times^HX:=(G\times X)/H$ where the quotient is taken with respect to the free $H$-action $h\cdot(g,x)=(gh^{-1},h\cdot x)$. 
In this case, $(G\times X)/H\simeq(G\times X)/\!/H$. The class of $(g,x)$ is denoted by $g*x$. The group $G$ naturally acts on $G\times^H X$ by left multiplication on the first component. 
It follows from the proof of  \cite[Lemma I.3]{luna} that if $Y\subset G\times^HX$ is  $G$-stable  and Zariski open, respectively closed, respectively locally closed, then  there exists  a $H$-stable  open, respectively closed, respectively locally closed subset $Y_X\subset X$ such that $Y=G\times^HY_X$. Also, there is a natural correspondence between $G$-orbits in $G\times^HX$ and $H$-orbits in $X$.

\begin{proposition}\label{prop:maffei}Let $r\in G$ be a semisimple element and let $M=C_G(r)$. There is a Zariski open neighbourhood $U$ of $r$ in $M$ such that:
\begin{enumerate}
\item $U$ is $\pi_{M}$-saturated; 
\item For any Jordan class $J_M$ of $M$ we have $J_M\cap U\neq \emptyset$ if and only if $r\in\overline{J_M}$;
\item The restriction $\gamma_U$ to $G\times^M U$ of the map $\gamma\colon G\times^M M\to G$ given by $\gamma(g*x)=gxg^{-1}$ is \'etale;
\item The image  $G\cdot U$ of $\gamma_U$ is a $\pi_G$-saturated open neighbourhood of $r$ in $G$.
\end{enumerate}
\end{proposition} 
\pf Observe that $G\cdot (1*r)=G*r$ and $\Oc_r^G$ are closed because $r$ is semisimple. By construction, the restriction of $\gamma$ to $G*r$ is injective. We claim that  $\gamma$ is \'etale at $1*r$, that is, the differential ${\textrm d}\gamma_{(1*r)}\colon T_{1*r}(G\times^{M}M)\to T_r G$ is bijective. 
We consider the map $\widetilde{\gamma}\colon G\times M\to G$ given by the conjugation action and the natural projection $p\colon G\times M\to G\times^M M$, so $\widetilde{\gamma}=\gamma\circ p$. 
For $m\in M$ the differential ${\textrm d}\widetilde{\gamma}_{(1,m)}\colon\gg\oplus\mathfrak{m}\to\gg$ at $(1,m)$ is given by $(x,y)\mapsto y-x+Ad(m^{-1})x$. For $g\in G$, let $L_g$ be left translation in $G$ by $g$. The induced map  identifies $\gg$ with $T_gG$ and $\mathfrak{m}=\cc_\gg(r)$ with $T_g M$. This way, ${\textrm d}\widetilde{\gamma}_{(g,m)}\colon\gg\oplus\mathfrak{m}\to\gg$ is given by $(x,y)\mapsto Ad(g)(y-x+Ad(m^{-1})x)$. Since $r$ is semisimple, $\gg={\rm Im}(Ad(r^{-1})-\id)\oplus{\rm Ker}(Ad(r^{-1})-\id)$ and 
${\rm Ker}(Ad(r^{-1})-\id)={\rm Ker}(\id- Ad(r))=\mathfrak{m}$ so ${\textrm d}\widetilde{\gamma}_{(1,r)}$ is onto, yielding surjectivity of  ${\textrm d}{\gamma}_{1*r}$.
For any pair $(g,m)\in G\times M$ the composition 
$$
\begin{CD}
G\times M@>{L_g\times L_m}>>G\times M@>p>>G\times^{M}M
\end{CD}
$$
yields an exact sequence
$$
\begin{CD}
0@>>>N_m@>>>\gg\oplus\mathfrak{m}@>{\textrm d}p>>T_{g*m}(G\times^{M}M)@>>>0
\end{CD}
$$
where $N_m=\{(x,x-Ad(m^{-1})(x)~|~x\in\mathfrak{m}\}$, so $\dim T_{g*m}(G\times^M M)=\dim \gg$ and injectivity of ${\textrm d}{\gamma}_{1*r}$ follows.
Therefore the hypotheses of \cite[Lemme fondamental, \S II.2]{luna} are satisfied for the map $\gamma\colon G\times^M M\to G$ and the point $1*r$ and there exists an \'etale slice; in particular, there exists a $\pi_M$-saturated Zariski open neighbourhood  $U'$ of $r$ in $M$ such that the restriction of $\gamma$ to $G\times^M U'\to G$ is \'etale with image a $\pi_G$-saturated open subset $V'=G\cdot U'$ of $G$.  

Consider the stratification on $M/\!/M$ with finitely-many closed strata of the form $\overline{J_M}/\!/M$, for $J_M$ a (semisimple) Jordan class in $M$, and let $J_M/\!/M$ denote the open stratum in $\overline{J_M}/\!/M$. 
Let $V$ be the union of all $J_M/\!/M$ containing the class of $r$ in their closure. It is open, because its complement is the closed set
\begin{equation*}\bigcup_{[r]\not\in \overline{J_M}/\!/M}J_M/\!/M=\bigcup_{[r]\not\in \overline{J_M}/\!/M}\overline{J_M}/\!/M.\end{equation*}
Then $U'':=\pi^{-1}_M(V)$ is a $\pi_M$-saturated open subset of $M$  containing $r$ and such that a Jordan class $J_M$ in $M$ meets $U''$ if and only if $r\in\overline{J_M}$. We take the $\pi_M$-saturated neighbourhood $U=U'\cap U''$. It satisfies condition 2. and the restriction of the \'etale map $\gamma$  to the open subset $G\times^M U$ is again \'etale and its image $G\cdot U$ is a $\pi_G$-saturated open neighbourhood of $r$ in $G$. 
\epf

\begin{remark}\label{rk:centraliser}{\rm With notation as above, since $\gamma_U$ is \'etale, for any $x\in U$ we have $\dim G\cdot x=\dim G\cdot(\gamma(1*x))=\dim(G*x)$, so $\dim C_G(x)=\dim G_{1*x}=\dim (C_G(x)\cap M)$. Hence, $C_M(x)^\circ=C_G(x)^\circ$. Since $U$ is $\pi_M$-saturated, if $x=su\in U$, then $s\in U$ and so $C_M(s)^\circ=C_G(s)$, see also  \cite[Remarque III.1.4]{luna}. }
\end{remark}

\begin{proposition}\label{prop:local}Let $J=J_G(\tau)$ for some $\tau\in \mathcal{T}$,  let $rv\in\overline{J}$ and let $M=C_G(r)$. Then 
\begin{align}\label{eq:se}
(\overline{J},\,rv)\seq\left(\bigcup_{i\in I_{J,rv}}r^{-1}\overline{J_{M,i}},\;v\right)
\end{align} where the $J_{M,i}$'s for $i\in I_{J,rv}$ are precisely the Jordan classes in $M$ contained in $J$ and containing $rv$ in their closure. 

If, in addition, $rv\in\overline{J}^{\reg}$, then
$rv\in\overline{J_{M,i}}^{\reg}$ for every $i\in I_{J,rv}$ and 
\begin{align}\label{eq:se-reg}
(\overline{J}^{\reg},\,rv)\seq\left(\bigcup_{i\in I_{J,rv}}r^{-1}\overline{J_{M,i}}^{\reg},\;v\right)\!.
\end{align}
\end{proposition}
\pf Let $\tau=(M', Z(M')^\circ s, \Oc)$. Since conjugation by $g\in G$ induces an isomorphism of pointed varieties $(\overline{J},rv)\simeq(\overline{J},g\cdot (rv))$, we may assume that $r\in Z(M')^\circ s$, so $M'\subset M$. We adopt notation from Proposition \ref{prop:maffei} and its proof, but with $\gamma_U$ viewed as a map $G\times^MU\to G\cdot U$. Let $\widetilde\gamma_U\colon G\times U\to G\cdot U$ be the restriction of $\widetilde{\gamma}$.  

We will first show that $(\overline{J},x)\seq(\overline{J}\cap U, x)$ for any $x\in \overline{J}\cap U$. Then, we will prove that $\overline{J}\cap U=\overline{J\cap U}^U$ and show that the irreducible components of $\overline{J\cap U}^U$ are the intersections of $U$ with the closures of those Jordan classes in $M$ that are contained in $J$ and contain $r$ in their closure. We will conclude the proof of \eqref{eq:se} by observing that, in order to study $(\overline{J}, x)$ we can neglect those irreducible components of $\overline{J\cap U}^U$ not containing $x$. A dimension argument will give \eqref{eq:se-reg}.   

We consider the following commutative diagram
$$
\begin{CD}
G\times U   @ >p >> G\times^MU@>{\gamma_U}>>G\cdot U\\
@AAA  @AAA@AA\iota A\\
G\times(\overline{J}\cap U)      @>>>
G\times^M(\overline{J}\cap U) @>>>\overline{J}\cap G\cdot U
\end{CD}
$$
  
Observe  that $\widetilde{\gamma}_U^{-1}(\overline{J}\cap G\cdot U)$ is a $G$-stable closed subset of $G\times U$, so it is of the form $G\times V$ for some $V$ closed in $U$. In turn, $V$ is the pre-image of $G\times V$ through the natural inclusion of $U$ into $G\times U$. Therefore $\widetilde{\gamma}_U^{-1}(\overline{J}\cap G\cdot U)=G\times(\overline{J}\cap U)$. 
This is exactly saying that  the bottom composition of arrows is obtained by pulling-back $\widetilde{\gamma}_U$ along  the closed embedding $\iota$.
Hence the bottom composition is also smooth and for any $x\in \overline{J}\cap U$
\begin{align}\label{eq:prime}&(\overline{J},x)\loc(\overline{J}\cap G\cdot U,x)\seq(G\times(\overline{J}\cap U), (1,x))
\seq(\overline{J}\cap U, x).\end{align} 
 We show that  $\overline{J}\cap U=\overline{J\cap U}^U$
by 
proving the equivalent statement $G\times^M(\overline{J}\cap U)=\overline{G\times^M(J\cap  U)}^{G\times^M U}$, i.e.,  $\gamma_U^{-1}(\overline{J}\cap G\cdot U)=\overline{\gamma_U^{-1}(J\cap G\cdot U)}^{G\times^M U}$.  By elementary topology we see that $\overline{J}\cap (G\cdot U)=\overline{J\cap G\cdot U}^{G\cdot U}$. Since $\gamma_U$ is continuous, surjective and open, 
$\gamma_U^{-1}(\overline{J\cap G\cdot U}^{G\cdot U})=\overline{\gamma_U^{-1}(J\cap G\cdot U)}^{G\times^M U}$ giving the desired equality. Thus \eqref{eq:prime} gives $(\overline{J},x)\seq(\overline{J\cap U}^U, x)$ for any $x\in \overline{J}\cap U$.

We describe now the irreducible components of $\overline{J\cap U}^U$. By base-change the restriction of $\gamma$ to $G\times^M(J\cap U)$ is  a $G$-equivariant \'etale map to $J\cap G\cdot U\subset G_{(d)}$ for some $d$. Hence all $G$-orbits in $G\times^M(J\cap U)$ have the same dimension. By Remark \ref{rk:centraliser}  we have  $J\cap U\subset M_{(d')}$ for some $d'$.  The equivalence \eqref{eq:prime} implies that $(J\cap U,x)\seq(J,x)$ for any $x\in U\cap J$ and $J$ is smooth, so the intersection $U\cap J$ is also smooth. Hence it is the union of its connected components and it is contained in the finite union of  those Jordan classes in $M_{(d')}$ containing $r$ in their closure. 
Let $J_M$ be a Jordan class in $M$ such that $J\cap U\cap J_M\neq\emptyset$. By construction of $U$, we have $r\in\overline{J_M}$.  It follows from Remark \ref{rk:centraliser} that if $x=tu\in J_M\cap U\cap J$, then $C_M(t)^\circ=C_G(t)$, hence $\dim Z(C_M(t)^\circ)^\circ=\dim Z(C_G(t)^\circ)^\circ=\dim Z(M')^\circ$.
The proof of \cite[Theorem 5.6 (e)]{gio-espo} applied to the case of (regular closures of) arbitrary Jordan classes shows that $\dim J_M=d'+\dim Z(M')^\circ$, so all Jordan classes of $M$ meeting $J\cap U$ have the same dimension.  The same argument also shows that  $(Z(C_M(t)^\circ)^\circ r)^{\reg} u=(Z(C_G(t))^\circ r)^{\reg}u$ and  so $J_M=M\cdot((Z(C_M(t)^\circ)^\circ r)^{\reg}u)\subset G\cdot((Z(C_G(t))^\circ r)^{\reg}u)=J$. Therefore, $J_M\subset J$. Conversely, if a Jordan class $J_M\subset M$ contains $r$ in its closure and is contained in $J$, then $\emptyset\neq J_M\cap U\subset J\cap U$. 

Let $I_{J,r}$ be the index set parametrising  the Jordan classes $J_{M,i}$ of $M$ such that $r\in\overline{J_{M,i}}$ and $J_{M,i}\subset J$. Then $J\cap U=\bigcup_{i\in I}J_{M,i}\cap U$, and the locally closed subsets $J_{M,i}\cap U$ are finitely many,  disjoint, irreducible and have all the same dimension. Hence, their closures are the irreducible components of $\overline{U\cap J}^U=\overline{J}\cap U$. 
Therefore,  for any $x\in \overline{U\cap J}^U$:
\begin{align*}(\overline{U\cap J}^U,x)&\simeq(\overline{\bigcup_{i\in I}J_{M,i}\cap U}^U,x)\simeq(\bigcup_{i\in I}\overline{J_{M,i}\cap U}^U,x)\\
&\simeq(\bigcup_{i\in I}\overline{J_{M,i}}^M\cap U,x)\loc(\bigcup_{i\in I}\overline{J_{M,i}}^M,x).\end{align*}
Let $I_{J,x}$ be the set of indices in $I_{J,r}$ such that $x\in\overline{J_{M,i}}$ and let $U_x$ be a Zariski open neighbourhood of $x$ in $M$ such that $U_x\cap\overline{ J_{M,i}}=\emptyset$ for any $i\in I_{J,r}\setminus I_{J,x}$. 
Then, 
\begin{align*}(\overline{J},x)\seq(\bigcup_{i\in I_{J,r}}\overline{J_{M,i}}^M\cap U_x,x)\loc(\bigcup_{i\in I_{J,x}}\overline{J_{M,i}}^M\cap U_x,x)\loc(\bigcup_{i\in I_{J,x}}\overline{J_{M,i}}^M,x).\end{align*} Taking $x=rv$ and translating by $r^{-1}$  gives \eqref{eq:se}. Observe that if $rv\in\overline{J}^{\reg}$, then $\Oc_{rv}^G\subset G_{(d)}$ and meets $U$. Since $\gamma_U$ is \'etale, $\Oc_{rv}^M\subset M_{(d')}$ so it lies in $\overline{J}_{M,i}^{\reg}$ for every $i\in I_{J,x}$. Since $\overline{J}^{\reg}$ is open in $\overline{J}$ and 
$\bigcup_{i\in I_{J,rv}}  \overline{J}_{M,i}^{\reg}$ is open in the union of equidimensional closures $\bigcup_{i\in I_{J,rv}}  \overline{J}_{M,i}$, equation \eqref{eq:se-reg} follows from \eqref{eq:se}.
\epf

\medskip

In order to provide an explicit parametrisaton of  the set $I_{J, rv}$ from Proposition \ref{prop:local} in terms of data depending on $J$ and $rv$, we introduce some notation.
Let $\tau=(M',Z(M')^\circ s,\Oc)\in{\mathcal T}$, let $rv\in\overline{J}\cap Z(M')^\circ sv$ and let $M=C_G(r)$. We set:
\begin{align*}
W_\tau&:={\rm Stab}_W(\tau)\\
W(\tau,r)&:=\{w\in W~|~r\in w\cdot (Z(M')^\circ s)\}.
\end{align*}
If $w\in W(\tau,r)$ then $w\cdot M'=C_G(w\cdot(Z(M')^\circ s))^\circ\subset M$ and $w\cdot M'$ is a Levi subgroup in $M$, \cite[Lemma 4.10]{gio-espo}. We consider then
\begin{align*}
W(\tau,rv)&:=\left\{w\in W(\tau,r)~|~\Oc_v^M\subset\overline{{\rm Ind}_{w\cdot M'}^M(w\cdot \Oc)}\right\}.
\end{align*}
The reader should be alerted that $W(\tau,r)$ and $W(\tau,rv)$ are  not  subgroups of $W$ in general.

Since $W_\tau\leq {\rm Stab}_W(Z(M')^\circ s)$, it acts on $W(\tau,r)$ from the right. It preserves $M'$ and $\Oc$, hence it acts also on $W(\tau,rv)$ from the right. \\
The group $W_r:=N_{M}(T)/T\leq W$ acts on  $W(\tau,r)$ and $W(\tau,rv)$ from the left. 

\begin{theorem}\label{thm:main}Let $J=J_G(\tau)$ for some $\tau=(M',Z(M')^\circ s,\Oc)\in\mathcal{T}$, let  $r\in Z(M')^\circ s$ and $M=C_G(r)$.  Then
\begin{align}\label{eq:Main}
(\overline{J},r)\seq\left(\bigcup_{w \in W_r\backslash W(\tau,r)/W_\tau} \overline{J_{M}(w\cdot\tau)},\,r\right).
\end{align}
If $rv\in\overline{J}$ then 
\begin{align}\label{eq:Main2}
(\overline{J},rv)\seq\left(\bigcup_{w \in W_r\backslash W(\tau,rv)/W_\tau} r^{-1}\overline{J_{M}(w\cdot\tau)},v\right).
\end{align}
If $rv\in\overline{J}^{\reg}$ then 
\begin{align}\label{eq:Main2-reg}
(\overline{J}^{\reg},rv)\seq\left(\bigcup_{w \in W_r\backslash W(\tau,rv)/W_\tau} r^{-1}\overline{J_{M}(w\cdot \tau)}^{\reg},v\right).
\end{align}
\end{theorem}
\pf We first consider the neighbourhood of $r$. By Proposition \ref{prop:local} it is enough to show that the right hand side of \eqref{eq:Main} involves precisely those Jordan classes in $M$ that 
\begin{enumerate}
\item are contained in $J$ and 
\item contain $r$  in their closure. 
\end{enumerate}
By condition 1, the latter are parametrised by $W_r$-orbits of triples of the form $w\cdot\tau$ for some $w\in W/W_\tau$. Condition 2 is equivalent to $r\in w\cdot (Z(M')^\circ s)$. Hence 
the elements $w$ must be taken in $W(\tau,r)/W_\tau$. This gives 
\eqref{eq:Main}. 

Let us now consider the neighbourhood of $rv$. In this case we need to prove that the classes occurring in the right hand side of \eqref{eq:Main2} are precisely those Jordan classes $J_M(M'',Z(M'')^\circ s',\Oc')$ in $M$ that 
\begin{enumerate}
\item are contained in $J$ and 
\item contain $rv$ in their closure, that is, contain $r$ in their closure and  satisfy $\Oc_v^M\subset\overline{{\rm Ind}_{M''}^M\Oc'}$. 
\end{enumerate}
They are parametrised by $W_r$-orbits of  triples  of the form $w\cdot\tau$, where $w$ must be taken in $W(\tau,rv)/W_\tau$, as one sees from condition 2.
This gives \eqref{eq:Main2}. Equation \eqref{eq:Main2-reg} follows from \eqref{eq:Main2} and \eqref{eq:se-reg}.
\epf

\begin{corollary}\label{cor:geometry-constant} Let $J=J_G(\tau)$ for some $\tau\in\mathcal{T}$ and let $rv,\,r'v'\in J'\subset\overline{J}$.
Then, $(\overline{J},rv)\seq(\overline{J},r'v')$. In other words, the geometry of $G$ and $\overline{J}$ is constant along  Jordan classes.
\end{corollary}
\pf Let $\tau=(M', Z(M')^\circ s, \Oc)$. Since $(\overline{J},x)\seq(\overline{J},g\cdot x)$ for any $g\in G$, we may assume that  $r\in Z(M')^\circ s$,  $C_G(r)=C_G(r')$,  $r'\in (Z(C_G(r))^\circ r)^{\reg}$ and $v'=v$ so $W_{r'}=W_r$. We set $M:=C_G(r)$. If $r\in w\cdot(Z(M')^\circ s)$ for some $w\in W$, then  $M\supset C_G(w\cdot(Z(M')^\circ s))=w\cdot M'$ whence $Z(M)^\circ \subset w\cdot Z(M')^\circ$, and therefore  $r'\in Z(M)^\circ r\subset w\cdot (Z(M')^\circ s)$. Hence, $W(\tau,r)=W(\tau,r')$ and so $W(\tau,rv)=W(\tau,r'v)$. 
The statement follows from \eqref{eq:Main2}  and left translation by $r'r^{-1}\in Z(M)^\circ$.
\epf

\begin{corollary}\label{cor:equivalent}Let $J=J_G(\tau)$, for $\tau=(M',Z(M')^\circ s,\Oc)\in\mathcal{T}$, let $rv\in \overline{J}\cap Z(M')^\circ sv$ and let $M=C_G(r)$. Then $\overline{J}$ is unibranch, respectively smooth, respectively normal, at $rv$ if and only if $\left|W_r\backslash W(\tau,rv)/W_\tau\right|=1$ and $r^{-1}\overline{J_{M}(\tau)}$ is so at $v$. 
\end{corollary}
\pf Let $U$ be  as in the proof of Proposition \ref{prop:local}. Then the irreducible components of $U\cap\overline{J}$ containing $rv$ are precisely the subsets 
$\overline{J_{M}(w\cdot \tau)}\cap U$ for $w\in W_r\backslash W(\tau,rv)/W_\tau.$
Hence, $\left|W_r\backslash W(\tau,rv)/W_\tau\right|=1$
 is a necessary condition for $\overline{J}$  being unibranch at $rv$, and a fortiori, normal, or smooth. In addition, if  $\left|W_r\backslash W(\tau,rv)/W_\tau\right|=1$,  then 
$(\overline{J},rv)\seq(r^{-1}\overline{J_{M}(\tau)},v)$ and we use the properties of smooth equivalence. 
\epf

The same argument gives the following statement.
\begin{corollary}\label{cor:CM}Let $J=J_G(\tau)$, for $\tau=(M',Z(M')^\circ s,\Oc)\in\mathcal{T}$, let $rv\in \overline{J}\cap Z(M)^\circ sv$ and let $M=C_G(r)$. Assume $\left|W_r\backslash W(\tau,rv)/W_\tau\right|=1$. Then $\overline{J}$ is Cohen-Macaulay at $rv$ if and only if $r^{-1}\overline{J_{M}(\tau)}$ is so at $v$. \hfill$\square$
\end{corollary}

The local study of the closure of a Jordan class $J=J_G(\tau)$ around $rv$ simplifies drastically when $\left|W_r\backslash W(\tau,r)/W_\tau\right|=1$ and therefore it is important to characterize when this is the case. The next corollary deals with this question under the assumption that $W_\tau={\rm Stab}_W(Z(M')^\circ s)$, which is always satisfied when $\Oc$ is characteristic, e.g., when $\Oc=1$ (semisimple Jordan classes) or when $\Oc$ is regular (regular Jordan classes). 

\begin{lemma}\label{lem:estimate}Let $J=J_G(\tau)$ for $\tau=(M', Z(M')^\circ s,\Oc)\in\mathcal{T}$ and let $r\in (Z(M')^\circ s)\cap\overline{J}$. Assume that $W_\tau={\rm Stab}_W(Z(M')^\circ s)$. Then $\left|W_r\backslash W(\tau,r)/W_\tau\right|=1$ if and only if $\overline{J}/\!/G$ is unibranch at the class $[r]$ of $r$.
\end{lemma}
\pf The isomorphism $G/\!/G\simeq T/W$ identifies $\overline{J}/\!/G$ with $W\cdot (Z(M')^\circ s)/W$, so we need to understand the neighbourhood of $W\cdot (Z(M')^\circ s)/W$ around $[r]$. By \cite[Anhang zu K. 7, Satz 21]{tedesco},  there is a $W_r$-stable analytic open neighbourhood $U$ of $r$ in $W\cdot(Z(M')^\circ s)$
such that $U/W_r$ identifies  with a neighbourhood of $[r]$ in $W\cdot(Z(M')^\circ s)/W$. We can choose $U$ so that it meets only the $W$-translates of $Z(M')^\circ s$ containing $r$. 
Therefore 
\begin{align*}
\left(W\cdot(Z(M')^\circ s)/W, [r]\right)&\loc\left(W\cdot(Z(M')^\circ s)\cap U/W_r, [r]\right)\\
&\loc\left(\bigcup_{w\in W(\tau,r)/W_\tau}w\cdot(Z(M')^\circ s)/W_r,[r]\right).\end{align*} 
Here, $W_r$ acts as usual from the left. Hence, $\overline{J}/\!/G$ is unibranch at $[r]$ if and only if $\left|W_r\backslash W(\tau,r)/W_r\right|=1$.
\epf

 \begin{remark}{\rm By construction $\left|W_r\backslash W(\tau,rv)/W_\tau\right|\leq \left|W_r\backslash W(\tau,r)/W_\tau\right|$ but the inequality may be strict: here is an example.
 Let $G={\rm SL}_{4}({\mathbb C})$,  $M=\langle T,\,X_{\pm\alpha_1}\rangle$, $\tau=(M, Z(M)^\circ, 1)$.  In this case $Z(M)=Z(M)^\circ$ and $W_\tau={\rm Stab}_W(Z(M))=\langle s_1,s_3\rangle$. 
Let $rv\in \overline{J_G(\tau)}$ with $C_G(r)=\langle T,\,X_{\pm\alpha_1},\,X_{\pm\alpha_3}\rangle$ and $v\in{\rm Ind}_{M}^{C_G(r)}(1)$. Then $W_r=W_\tau$ and $v$ is trivial in the component corresponding to $\alpha_1$ and regular in the component corresponding to $\alpha_3$. If $w\in W$ satisfies $r\in w \cdot Z(M) \neq Z(M)$, then $w\in w_0 W_\tau$. Since $w_0\not\in W_\tau$
we have $\left|W_r\backslash W(\tau,r)/W_\tau\right|=2$. However, if $w$ is as above,  then $rv\not\in\overline{J_{C_G(r)}(w \cdot\tau)}$. Indeed, if $rv'\in \overline{J_{C_G(r)}(w\cdot\tau)}$, then $v'\in  \overline{{\rm Ind}_{w \cdot M}^{C_G(r)}(1)}$ which does not contain ${\rm Ind}_{M}^{C_G(r)}(1)$.  Hence 
$\left|W_r\backslash W(\tau,rv)/W_\tau\right|=1$.}
 \end{remark}
 
\section{From unipotent elements to nilpotent elements}\label{sec:toLie}

In this section $G$ is an arbitrary complex connected reductive group. We compare the local geometry of Jordan classes containing a unipotent element with the local geometry of Jordan classes in $\gg$. 
Our main tool will be the exponential map, as introduced in Section \ref{sec:notation}.
It is well-known that it is an analytic local isomorphism around a nilpotent element, see for instance \cite[Chapter I, Theorem 3.5]{GV}. 

For our purposes we will need to describe more explicitly the loci on which $\exp_{sc}$ and $\exp$ are local analytic isomorphisms as open neighbourhoods of the nilpotent cone $\NN$. We will use the convention that, unless otherwise stated,  a letter in gothic character will denote the Lie algebra of the group denoted by the same letter in capital Latin character.

If $G$ is simple, we consider the coroot lattice  $Q^\vee={\rm Ker}(\exp_{sc})$; its real span $\hh_{\mathbb R}={\mathbb R}Q^\vee$, so $\hh=\hh_{\mathbb R}\otimes_{\mathbb R}{\mathbb C}$; the fundamental alcove $\mathcal{A}:=\{h\in\hh_{\mathbb R}~|~0\leq\alpha(h)\leq1,\;\;\forall\alpha\in\widetilde{\Delta}\}$ and the affine hyperplanes $H_{\alpha, l}:=\{h\in \hh~|~\alpha(h)=l\}$ for $l\in{\mathbb Z}$. Let $A$ be the interior of $W\cdot\mathcal{A}+i\hh_{\mathbb R}$ in $\hh$ and let 
\begin{equation}\label{eq:U}
U_{sc}:=\pi_\gg^{-1}\pi_\gg A.
\end{equation} 
It is $\pi_\gg$-saturated by construction and open by Chevalley's restriction theorem. 
\begin{lemma}\label{lem:Lie}
There exist a $\pi_{\gg}$-saturated analytic open neighbourhood $U$ of $\NN$ in $\gg$ and a $\pi_{G}$-saturated analytic open neighbourhood $V$ of $\UU$ in $G$ such that 
the restriction of $\exp$ to $U$  is an analytic isomorphism $\exp_U\colon U\to V$. If $G=G_{sc}$ and it is simple, then one can take $U=U_{sc}$ as in \eqref{eq:U} and $V=\exp_{sc}(U_{sc})$. 
\end{lemma}
\pf If $G$ is a torus, then $\gg$ is abelian, $\exp$ is a local analytic isomorphism and $\NN$ and $\UU$ are trivial, so there is nothing to prove. If $G$ is a direct product, then it is enough to prove the statement for each factor. 
Let $\pi\colon G_{sc}\to G$ be a central isogeny and let $\overline{\pi}\colon G_{sc}/\!/G_{sc}\to G/\!/G$ be the induced map and assume there exist a $\pi_{\gg}$-saturated analytic open neighbourhood $\widetilde U$ of $\NN$ in $\gg$ such that the restriction of $\exp_{sc}$ to $\widetilde{U}$  is an analytic isomorphism. Let $A'$ be an open neighbourhood of the class $[1]$ in  $G_{sc}/\!/G_{sc}$ such that if  $kA'\cap A'\neq\emptyset$ for some $k\in {\rm Ker}\,\overline{\pi}$, then $k=1$. Let $\widetilde{A}=\pi_{G_{sc}}^{-1}A'$. Then $\widetilde{V}:=\widetilde{A}\cap\exp_{sc}\widetilde{U}$ is a $\pi_{G_{sc}}$-saturated open neighbourhood of $\UU$ in $G_{sc}$ and $U:=\widetilde{U}\cap \exp_{sc}^{-1}(\widetilde{V})$ and $V:=\exp(U)=\pi(\widetilde{V})$ are the sought neighbourhoods for $\gg$ and $G$.

Since there is always a central isogeny $\pi\colon Z(G)^\circ\times [G,\,G]_{sc}\to G$ and $[G,\,G]_{sc}$ is a direct product of simple, simply connected factors, it remains to prove the statement for $G$  simple and simply connected. In this case we show that $\exp_{sc}$ is an analytic isomorphism on $U_{sc}$. 
The main result in \cite{nono,nono2}  (see also \cite[Chapter I, Theorem 3.5]{GV}), states that the exponential map is a local analytic isomorphism at $x=x_s+x_n$ if and only if the eigenvalues of $\ad(x)$ do not meet $\mathbb{Z}\setminus\{0\}$. These eigenvalues coincide with those of $\ad(x_s)$, so the condition is verified if and only if, up to $G$-action, $x_s$ lies in $\hh\setminus\bigcup_{l\in {\mathbb Z}^\times,\alpha\in\Phi^+}H_{\alpha, l}$. As $A$ is contained in this set, $\exp_{sc}$ is a local analytic isomorphism on $U_{sc}$. 
Let $V=\exp_{sc}(U_{sc})$. We prove now that the restriction $\exp_{U}\colon U\to V$ is an analytic isomorphism , i.e., that $\exp_{sc}$ is injective on $U_{sc}$. If $\exp_U(x_s+x_n)=\exp_U(y_s+y_n)$, then $x_n=y_n$ because $\exp_{sc}$ is an isomorphism on $\NN$ and by $G$-equivariance we may assume that $\exp_U(x_s)=\exp_U(y_s)\in T$, so $x_s,\,y_s\in A$. Two elements in $A$ cannot differ by an element in $Q^\vee$ because ${\mathcal A}$ is a fundamental domain for the $Q^\vee\rtimes W$-action on $\hh_{\mathbb R}$, \cite[Theorem 4.8]{hu-rgcg} and $Q^\vee$ does not change the imaginary components of elements in $\hh$. Thus, $x_s=y_s$ and $x=y$. The properties of $V$ follow from those of $\exp_{sc}$ and of $U_{sc}$.
\epf

We describe now compatibility of the Jordan stratifications induced on $U$ and $V$ when $U$ and $V$ are as above. 

\begin{theorem}\label{thm:stratification}
Let $U$ be a $\pi_{\mathfrak g}$-saturated analytic open neighbourhood of $\NN$ in $\gg$ and  let $V$ be a $\pi_{G}$-saturated analytic open neighbourhood of $\UU$ in $G$ such that the restriction of $\exp$ to $U$ is a $G$-equivariant analytic isomorphism $\exp_U\colon U\to V$. Then, $\exp_U$  identifies the stratification on $U$ induced by the Jordan one in $\gg$ with the stratification on $V$ induced by the Jordan one in $G$, preserving dimensions, closure orderings, orbit dimensions. More precisely, for $\tau=(M,Z(M)^\circ s,\Oc)\in\mathcal{T}$ we have $J_G(\tau)\cap V\neq\emptyset$ if and only if $M$ is a Levi subgroup of $G$ and $Z(M)^\circ s=Z(M)^\circ$ and if this is the case, then 
\begin{align*}J_G(M,Z(M)^\circ,\Oc)\cap V=\exp(\jj_\gg(\mathfrak{m},\mathfrak{O})\cap U).\end{align*}  
where $\exp\mathfrak{O}=\Oc$.
\end{theorem}
\pf We keep notation from the proof of Lemma \ref{lem:Lie}. Let $\jj=\jj_\gg(\ll,\mathfrak{O})$ be a Jordan class in $\gg$. Then $\overline{\jj}\cap {\mathcal N}\neq\emptyset$ so $\jj\cap U\neq\emptyset$. By $\pi_{\gg}$-saturation of $U$ we have
\begin{align*}U\cap(\zz(\ll)^{\reg}+\mathfrak{O})=U\cap\zz(\ll)^{\reg}+\mathfrak{O}.
\end{align*} 
If $x=x_s+x_n\in \zz(\ll)^{\reg}\cap U+\mathfrak{O}$, then  $\exp$ is a local diffeomorphism at $x$ by \cite{nono2}, rephrased in  \cite[Chapter I, Theorem 3.5]{GV}, so we have $\ll=\cc_\gg(x_s)=\cc_\gg(\exp(x_s))$ and $L=C_G(\exp(x_s))^\circ$ is a Levi subgroup of $G$. Setting  $\Oc=\exp \mathfrak{O}$ we have 
$\exp(\zz(\ll)^{\reg}\cap U+\mathfrak{O})\in V\cap Z(L)\Oc$.  Observe that $\zz(\ll)^{\reg}$ is obtained removing finitely many vector spaces of real codimension at least $2$ from a (complex) vector space, so it is connected in the analytic topology. Therefore $U\cap\zz(\ll)^{\reg}$, $U\cap\zz(\ll)^{\reg}+\mathfrak{O}$ 
and 
$\jj\cap U=Ad(G)(U\cap\zz(\ll)^{\reg}+\mathfrak{O})$ are also connected. By continuity, $\exp(U\cap(\zz(\ll)^{\reg}+\mathfrak{O}))$ and $\exp(U\cap\jj)$ are also connected in the analytic topology.
Thus, $\exp(\zz(\ll)^{\reg}\cap U+\mathfrak{O})\in V\cap (Z(L)^\circ s)^{\reg}\Oc$ for some $s\in Z(L)$ and $\exp(\jj\cap U)\subset J_G(L,Z(L)^\circ s,\Oc)\cap V$. Observe also that $0\in \overline{\jj}\cap U$ so 
 $1\in \overline{J(L,Z(L)^\circ s,\Oc)}\cap V$. By the discussion following \eqref{eq:closure-group} this implies that $Z(L)^\circ s=Z(L)^\circ$. 
 
Conversely, let $J$ be a Jordan class in $G$ such that $V\cap J\neq\emptyset$ and let $su\in V\cap J$, with $M=C_G(s)^\circ$. By $\pi_G$-saturation of $V$ we have
$(Z(M)^\circ s)^{\reg}\Oc_u^M\cap V=((Z(M)^\circ s)^{\reg}\cap V)\Oc_u^M$. For any $r\in (Z(M)^\circ s)^{\reg}\cap V$ we have $r=\exp(x_r)$ for some $x_r\in U$ and
 ${\rm Lie}(C_G(r)^\circ)=\mathfrak{m}=\cc_\gg(x_s)$. Therefore for any $rv\in (Z(M)^\circ s)^{\reg}\Oc_u^M)\cap V$ we have $rv\in \exp(U\cap\jj_\gg(\mathfrak{m},\mathfrak{O}_{\exp^{-1}u}^M))\subseteq J_G(M, Z(M)^\circ,\Oc_u^M)$ so
 $Z(M)^\circ s=Z(M)^\circ$ and  $\exp(U\cap\jj_\gg(\mathfrak{m},\mathfrak{O}_{\exp^{-1}u}^M))=V\cap J_G(M, Z(M)^\circ,\Oc_u^M)$.

Finally, $\exp_U$ is a $G$-equivariant analytic isomorphism, hence it preserves orbit dimensions, closure orderings, and dimensions.   
\epf

\begin{corollary}\label{cor:smooth_equivalence}Let $J=J_G(M,Z(M)^\circ,\Oc)$ with $v\in \overline{J}\cap\mathcal{U}$, let $\exp\mathfrak{O}=\Oc$ and $\exp x_n=v$. Then,  
\begin{align*}(\overline{J},v)&\loc(\overline{\jj_\gg(\mathfrak{m}, \mathfrak{O})},\,x_n)\\
(\overline{J}^{\reg},v)&\loc(\overline{\jj_\gg(\mathfrak{m}, \mathfrak{O} )}^{\reg},\,x_n).\end{align*} 
\end{corollary}
\pf Let $U$ and $V$ be neighbourhoods of $\mathcal N$ and $\mathcal U$, respectively, as in Lemma \ref{lem:Lie}, Theorem \ref{thm:stratification}. Then $v\in \overline{J}\cap V$
and $\exp_U$ is an analytic isomorphism mapping $\overline{\jj_\gg(\mathfrak{m},  \mathfrak{O})}\cap U$ to $\overline{J}\cap V$. 
\epf

\begin{remark}
{\rm \begin{enumerate}
\item The set of points $x$ in $\gg$ such that $\exp$ is a local analytic isomorphism at $x$ is not a union of Jordan classes in general.
For instance $s={\rm diag}(i , -i)$ and $s'={\rm diag}(1, -1)$  lie in the same Jordan class in $\mathfrak{sl}_2(\CC)$, and the condition on the eigenvalues in \cite{nono2} holds for $s$ but not for $s'$.
\item The image of $\exp$ is a union of Jordan classes in $G$. Indeed, $g=ru\in \exp \gg$ if and only if $r\in C_G(u)^\circ$, by \cite{Dj}. This condition is clearly $G$-stable, so it is enough to show that $r\in C_G(u)^\circ$ implies  $Z(C_G(r)^\circ)^\circ r\subset C_G(u)^\circ$. Now, $u\in C_G(r)^\circ$, so $Z(C_G(r)^\circ)\subset C_G(u)$.
Since $Z(C_G(r))^\circ r$ is connected and contains $r$, we have the desired inclusion. 
\end{enumerate}}
\end{remark}

\section{Applications}

In this Section $G$ is semisimple and,  simply connected and we  apply the results from Sections \ref{sec:to_uni} and \ref{sec:toLie} to deduce geometric properties of closures of regular Jordan classes, sheets and Lusztig strata.

\subsection{Closures of regular Jordan classes in $G$}

We recall that a Jordan class $J=J_G(M, Z(M)^\circ s, \Oc)$  in $G$ is called regular if $J\subset G^{\reg}$, i.e.,  if $\Oc=\Oc_{\reg}$, the regular unipotent class in $M$. 

\begin{theorem}\label{thm:equivalence} Let $J$ be a regular Jordan class in $G$. Then the following statements are equivalent:
\begin{enumerate}
 \item $\overline{J}$ is normal and Cohen-Macaulay. 
\item $\overline{J}$ is normal.
\item $\overline{J}/\!/G$ is normal.
\item $\overline{J}/\!/G$ is smooth.
\end{enumerate}
\end{theorem}
\pf Clearly 1$\Rightarrow$ 2 $\Rightarrow$ 3, see \cite[Paragraph 0.2]{mum-git} for the second implication. Also, 3 $\Leftrightarrow$ 4 by \cite[Corollary 8.1]{gio-espo-normal}.  We show that 3 $\Rightarrow$ 1. Let $J=J_G(\tau)$ for $\tau=(M', Z(M')^\circ s, \Oc_{\reg}^{M'})\in\mathcal{T}$. Recall that $\overline{J}/\!/G=\overline{J_G(M',Z(M')^\circ s,1)}/\!/G$. Let us assume $\overline{J}/\!/G$ is normal. Then it is everywhere unibranch and since the regular unipotent class is characteristic, Lemma \ref{lem:estimate} 
 gives $\left|W_r\backslash W(\tau,rv)/W_\tau\right|=1$ for all points $rv\in\overline{J}$. Since the locus where $\overline{J}$ is not normal (not Cohen-Macaulay, respectively) is closed, \cite[\href{https://stacks.math.columbia.edu/tag/00RD}{Tag 00RD}]{stacks-project} and the geometry of $\overline{J}$ is constant along Jordan classes  by Corollary \ref{cor:geometry-constant},  it is enough to check the desired properties of $\overline{J}$ at points in closed Jordan classes in $\overline{J}$. These are the  Jordan classes $J_G(M, Z(M)^\circ r,1)\subset\overline{J}$ with $M$ semisimple, i.e.,  isolated semisimple conjugacy classes in $G$,  see \S \ref{sec:Jordan_prelim}. Let $\Oc_r^G$ be such a class, with $M=C_G(r)$. 
By Corollaries \ref{cor:equivalent}, \ref{cor:CM} and \ref{cor:smooth_equivalence}, $\overline{J}$ is normal and Cohen-Macaulay at $r$ if and only if 
$\overline{\mathfrak{J}_{\mathfrak{m}}(\mathfrak{m}',\mathfrak{O}_{\reg}^{M'})}$ is so. 
By \cite[Theorem B]{richardson}, this happens if and only if  ${\rm Stab}_{W_r}(\mathfrak{z}(\mathfrak{m'}))$ acts on $\mathfrak{z}(\mathfrak{m'})$ as a reflection group and 
$\overline{\mathfrak{J}_{\mathfrak{m}}(\mathfrak{m}',0)}/\!/M$ is normal. The first condition is ensured by  \cite[Proposition 2.5, Lemma 8.3 (i)]{gio-espo-normal} applied to $\overline{\mathfrak{J}_{\mathfrak{m}}(\mathfrak{m}', 0)}/\!/M$. The second condition is ensured by \cite[Theorem 4.9]{gio-espo-normal} applied to $\overline{J}/\!/G$.
\epf

\begin{remark}
{\rm The fact that normality of $\overline{\mathfrak{J}_\mathfrak{m}(\mathfrak{m}',0)}/\!/M$ implies that ${\rm Stab}_{W_r}(\mathfrak{z}(\mathfrak{m}'))$ acts on $\mathfrak{z}(\mathfrak{m}')$ as a reflection group can also be deduced from the proof of \cite[Theorem 3.1]{broer} or from the main result in \cite{DR}.}
\end{remark}

\begin{remark}\label{cor:lista2}{\rm Let $G$ be simple. The Jordan classes in $G$ satisfying condition $3$ from Theorem \ref{thm:equivalence} are classified  in \cite[Theorem 8.7]{gio-espo-normal}. Therefore, the closure of a regular Jordan class  $J=J_G(M,Z(M)^\circ s,\Oc_{\reg}^M)$  is smooth if and only if $M$ is either $T$, or semisimple, or of the form $G_\Pi$ where  $\emptyset\subsetneq\Pi\subsetneq\widetilde{\Delta}$ is one of the subsets occurring in the mentioned classification.} \end{remark}

\subsection{Sheets}

In this Subsection we apply the local description to the case of sheets, i.e,  the regular closures of Jordan classes $J=J_G(M',Z(M')^\circ s,\Oc)$ with $\Oc$  rigid in $M'$.  We will apply repeatedly the following argument.
\begin{remark}\label{rk:reduction}{\rm Let $S=\overline{J}^{\reg}$, with $J=J_G(M,Z(M)^\circ s,\Oc)$ be a sheet in $G$.
\begin{enumerate}
\item The  locus where $S$ is not smooth, respectively normal,  is closed. Thus, by  Corollary \ref{cor:geometry-constant} it is enough to check smoothness or normality of $S$ at a point in each closed Jordan class  in $S$. These are Jordan classes of triples $(M', Z(M')^\circ s',\Oc')$ with $M'$ semisimple and are precisely the conjugacy classes of isolated elements contained in $S$, see \S \ref{sec:Jordan_prelim}. 
\item The conjugacy class  $w\cdot\Oc$ is rigid in $w\cdot M$ for any $w\in W$ and therefore \eqref{eq:Main2-reg} implies that $S$ in the neighbourhood of an isolated point $rv$ is smoothly equivalent to a union of sheets in the semisimple group $C_G(r)$ in the neighbourhood of $v$.  
\item As $\exp$ is compatible with induction, it maps rigid nilpotent orbits in $\gg$ to rigid unipotent conjugacy classes in $G$. Hence, it identifies a neighbourhood of $v$ in a sheet in $C_G(r)$ with a neighbourhood of a nilpotent element in a sheet of $\cc_\gg(r)$.
\end{enumerate}}
\end{remark}

\begin{theorem}\label{thm:sheet-smooth}Let $\Phi$ be classical and let $S=\overline{J(\tau)}^{\reg}$   be a sheet in  $G$. Then $S$ is smooth if and only if it is normal if and only if it is unibranch.
\end{theorem}
\pf One direction is immediate. Assume  $S$ is unibranch: we prove that it is smooth. Let $\tau=(M, Z(M)^\circ s,\Oc)$ and $\Oc=\exp\mathfrak{O}$. By Corollary \ref{cor:equivalent}  we have  $\left|W_r\backslash W(\tau,rv)/W_\tau\right|=1$ for any point $rv\in S$. Hence \eqref{eq:Main2-reg} and Corollary \ref{cor:smooth_equivalence} imply that $S$ is smooth at $rv$ if and only if 
$\overline{\mathfrak{J}_{\cc_\gg(r)}(\mathfrak{m},\mathfrak{O})}^{\reg}$ is smooth. By Remark \ref{rk:reduction} part 1, it suffices to prove smoothness of $S$ at isolated classes. In this case
$\cc_\gg(r)$ is semisimple and classical because its Dynkin diagram is a sub-diagram of the extended Dynkin diagram of $\gg$. In addition, $\overline{\mathfrak{J}_{\cc_\gg(r)}(\mathfrak{m},\mathfrak{O})}^{\reg}$ is a sheet in $\cc_\gg(r)$ by Remark \ref{rk:reduction}, part 2. Since all sheets in classical Lie algebras are smooth \cite{bongartz,peterson,imhof}, we have the statement.
\epf

\begin{theorem}\label{thm:sheet-smooth-cod-one}Let $\Phi$ be classical and  irreducible and let $S=\overline{J(\tau)}^{\reg}$  be a sheet in  $G$.
If $\overline{J(\tau)}/\!/G$ is normal in codimension $1$, then $S$ is smooth. 
 \end{theorem}
\pf By Theorem \ref{thm:sheet-smooth} it is enough to show that $S$ is unibranch at every isolated $rv\in S$. Let $\tau=(M,Z(M)^\circ s,\Oc)$. If $\overline{J(\tau)}/\!/G$ is normal in codimension $1$, then it is unibranch by \cite[Lemma 8.2, Lemma 8.3]{gio-espo-normal}.
By \cite[Lemma 3.3]{gio-espo-joseph} if $G$ is simple and simply connected and $\overline{J(\tau)}^{\reg}$ is a sheet we always have $W_\tau={\rm Stab}_W(Z(M)^\circ s)$, so Lemma \ref{lem:estimate} applies. 
\epf

\begin{remark}\label{cor:lista-lisce}{\rm  Let $\Phi$ be classical and  irreducible and let $\tau = (M, Z(M)^\circ s, \Oc) \in \mathcal{T}$. Now, $\overline{J(\tau)}/\!/G$ is normal in codimension $1$ if and only if $M$ is either $T$, semisimple, or equals $G_\Pi$ for some $\Pi \subset\widetilde{\Delta}$ occurring in the classification in \cite[Proposition 8.6]{gio-espo-normal}. If this is the case and $\Oc$ is rigid, then the sheet $S=\overline{J_G(M,Z(M)^\circ s,\Oc)}^{\reg}$  is smooth
by Theorem \ref{thm:sheet-smooth-cod-one}.}
\end{remark}


\begin{corollary}\label{cor:lista-lisce-ecc}
Let $G$ be simple with $\Phi$ exceptional. Let $M$ be either semisimple, $T$, or  $G_\Pi$ for   $\emptyset\subsetneq\Pi\subsetneq\widetilde{\Delta}$ of the following type:
\begin{itemize}
\item[$E_6$:] $A_5$, $D_4$, $4A_1$, $2A_2$,
\item[$E_7$:] $E_6$,  $D_6$, $D_4 + 2A_1$, $3A_2$, $2A_3$, $A_3 + 3A_1$, $D_4 + A_1$, $5A_1$, $\{\alpha_0,\alpha_1,\alpha_2,\alpha_3,\alpha_4\}$, $\{\alpha_2,\alpha_4,\alpha_5,\alpha_6,\alpha_7\}$,
 $D_4$, $\{ \alpha_0,\alpha_2, \alpha_3\}$, $\{ \alpha_2,\alpha_5, \alpha_7 \}$,  $\{ \alpha_0,\alpha_3, \alpha_5, \alpha_7 \}$.
\item[$E_8$:] $\widetilde{\Delta} \setminus \{\alpha_1, \alpha_3\}$, $\widetilde{\Delta} \setminus \{\alpha_1, \alpha_3, \alpha_6\}$, $\widetilde{\Delta} \setminus \{\alpha_4, \alpha_6, \alpha_8\}$, $\{\alpha_2, \alpha_5, \alpha_7, \alpha_0\}$, 
$D_7$, $E_7$,  $D_6 + A_1$, $2A_3 + A_1$, $3A_2 + A_1$, $D_5 + 2A_1$, $D_4 + A_3$, $D_6$, $E_6$, 
$D_4 + 2A_1$, $3A_2$, $D_4$,
\item[$F_4$:] $A_3$, $A_1 + B_2$, $2A_1 + \widetilde{A}_1$, $B_3$, $C_3$, $2A_1$, $\widetilde{A}_2$, $B_2$,
\item[$G_2$:] $\widetilde{A}_1$,
\end{itemize}
and let $\tau=(M,Z(M)^\circ s,\Oc)\in{\mathcal T}$ with $\Oc$ rigid in $M$. Then $S=\overline{J_G(\tau)}^{\reg}$  is smooth if and only if either $M=T$ or $M$ is semisimple or the pair $(\Pi,\Oc)$ is different from:
\begin{itemize}
\item[$E_7$:] $(D_6, [2^4,1^4])$,
\item[$E_8$:] $(E_7, 2A_1)$, $(E_7, (A_1 + A_3)a)$, $(D_6 + A_1, [2^4,1^4] + [1^2])$ and $(D_6, [2^4,1^4])$,
\item[$F_4$:] $(B_2,[1^5])$.
\end{itemize}
\end{corollary}
\pf If $M=T$ or $M$ is semisimple, then $S=G^{\reg}$ or a single conjugacy class and there is nothing to prove.  Let $M=G_{\Pi}$ with $\Pi$ from the above list. We apply Remark \ref{rk:reduction} and we look at $S$ in the neighbourhood of isolated elements $rv$. For all $\Pi$ the quotient  $\overline{S}/\!/G$ is normal in codimension 1, \cite[Proposition 8.6]{gio-espo-normal}, hence it is unibranch. In addition, \cite[Lemma 3.3]{gio-espo-joseph} ensures that $W_\tau={\rm Stab}_W(Z(M)^\circ s)$ for any choice of $Z(M)^\circ s$. By Lemma \ref{lem:estimate}, Corollary \ref{cor:equivalent} and \cite{imhof} the problem is reduced to showing that $\overline{\jj_{\mathfrak{m}'}(\mathfrak{m}, \mathfrak{O})}^{\reg}$ is smooth for $\Oc=\exp\mathfrak{O}$ and any ${\mathfrak{m}'}=\cc_\gg(r)$ semisimple exceptional containing  $\mathfrak{m}$.
Such Lie subalgebras are conjugate to $\gg_{\Pi'}$ for some $\Pi'\subset\widetilde{\Delta}$ with $|\Pi'|=|\Delta|$ and $\mathfrak{m}$ is $W_{\Pi'}$-conjugate to a standard Levi subalgebra therein, \cite[Lemma 4.9]{gio-espo}. 
However, normality in codimension 1 of $\overline{J}/\!/G$  is equivalent to the condition $\{w\Pi\subset\Phi~|~ w\in W,\,w\Pi\subset\widetilde{\Delta}\}=\{\Pi\}$. Therefore we are left to verify smoothness of the sheets  
$\overline{\mathfrak{J}_{\gg_{\Pi'}}(\gg_\Pi,\mathfrak{O})}^{\reg}$ for all exceptional $\Pi'\supset\Pi$ with $|\Pi'|=|\Delta|$. This is done by using the list in \cite[\S 4]{bulois} of smooth and singular sheets in simple exceptional Lie algebras on each simple component of  $\gg_{\Pi'}$.
\epf

\subsection{Sheets and Lusztig strata in ${\rm SL}_n({\mathbb C})$}

The case in which $G={\rm SL}_n({\mathbb C)}$ is particularly simple and we retrieve information on all its sheets and, as a consequence, on all Lusztig strata as defined in \cite[\S2]{lustrata}, see also \cite[\S3.2,3.3]{lustrata}.
\begin{proposition}\label{prop:sln}
Every sheet and Lusztig stratum in ${\rm SL}_n({\mathbb C)}$ 
is smooth.
\end{proposition}
\pf Let $S$ be a sheet in $G={\rm SL}_n({\mathbb C)}$. 
By Remark \ref{rk:reduction} 1., it suffices to prove smoothness at its isolated classes. These are all of the form $zv$ with $z$ central and $v$ unipotent, hence
$(S,zv)\seq(z^{-1}S,v)\seq(\mathfrak{S},\exp^{-1}v)$ where $\mathfrak{S}$ is a sheet in $\mathfrak{sl}_n$ by Corollary \ref{cor:smooth_equivalence} and Remark \ref{rk:reduction} 3.  All sheets in $\mathfrak{sl}_n$ are smooth by \cite{bongartz},  \cite{peterson}. Hence $S$ is smooth. 

We turn now to Lusztig strata. It follows from \cite[\S 2]{gio-MR} that their irreducible components are sheets in $G$. In the present case  strata are of the form $X_S=\bigcup_{k\in Z(G)}kS$ for $S=\overline{J_G(M,Z(M)^\circ s,1)}^{\reg}$ a given sheet.  We claim that $kS\cap k'S\neq\emptyset$ for some $k,k'\in Z(G)$ implies $kS=k'S$. 
Indeed,  a non-empty intersection of sheets always contains an isolated class \cite[Proposition 3.4]{gio-MR},  i.e, a class of the form $k''\Oc_v^G$ for $k''\in Z(G)$ and $v\in\mathcal{U}$. 
Observe that $kS=\overline{J_G(M,Z(M)^\circ ks,1)}^{\reg}$, for any $k\in Z(G)$. 
Formula \eqref{eq:reg-closure-group} gives $k''\in Z(M)^\circ ks\cap Z(M)^\circ k's$, i.e., $Z(M)^\circ ks=Z(M)^\circ k's$ and $kS=k'S$. Hence sheets in $G$ are connected components of strata so the latter are also smooth.
\epf

\section{Acknowledgements}
We are indebted to Andrea Maffei for interesting conversations and for the proof of Proposition \ref{prop:maffei}, and to Michael Bulois for communicating  the complete list of smooth sheets in Lie algebras of exceptional types, \cite{bulois}. We thank the referee for careful reading of the manuscript. This research was partially supported by DOR1898721/18,  DOR1717189/17  and BIRD179758/17 funded by the University of Padova.

\end{document}